\documentclass{article}
\usepackage{definitions}
\usepackage{graphicx} 
\usepackage{todonotes}

\newcommand{\mygmv}{\operatorname{\texttt{GMV}}}

\title{FlashFolio: A GPU-Accelerated Solver for Portfolio Optimization}
\author{Yilun Jiang\thanks{Point72 Asset Management (sjtujyl@gmail.com)} \and Haihao Lu\thanks{MIT, Sloan School of Management (haihao@mit.edu).} \and Zedong Peng\thanks{MIT, Sloan School of Management (zdpeng@mit.edu).} \and Jinwen Yang\thanks{University of Chicago, Department of Statistics (jinweny@uchicago.edu).}}
\date{}

\begin{document}

\maketitle

\begin{abstract}
We present FlashFolio, a GPU-accelerated solver for single-period and multi-period portfolio optimization with factor-based risk modeling, bid--offer spread costs, and nonlinear market impact. These models are widely used in portfolio construction and optimal execution, but become computationally challenging at large scale, especially in the multi-period setting. We benchmark FlashFolio against MOSEK on instances constructed from realistic market inputs. FlashFolio delivers consistent runtime improvements, achieving speedups of up to $12.9\times$ in the single-period setting and $48\times$ in the multi-period setting, while also exhibiting stronger robustness on challenging multi-period instances. Our results show that GPU-based optimization can help improve the practicality of large-scale portfolio optimization.

\end{abstract}

\section{Introduction}

Portfolio optimization is a central problem in operations research and financial engineering, concerned with allocating capital across a universe of assets to balance expected return, risk, and trading costs. Since the seminal mean-variance framework~\cite{markowitz-mv}, portfolio optimization has evolved into a flexible modeling paradigm that incorporates realistic features such as multifactor risk models, transaction costs, and portfolio constraints. These models play a critical role in a wide range of applications, including asset management and algorithmic trading. In modern financial markets, the increasing scale of tradable assets and the need for rapid decision-making have made efficient portfolio optimization an essential component of real-time trading and risk management systems.

A standard approach to portfolio optimization formulates the problem as a convex program, often a quadratic or conic optimization problem, where expected returns enter linearly, risk is modeled via a covariance matrix~\cite{FAMA19933,msciBarra2010}, and trading costs are captured through linear and nonlinear penalty terms~\cite{almgren2005direct, gatheral2010no}. Such formulations are attractive due to their interpretability, flexibility, and strong theoretical guarantees. In particular, multi-period extensions of these models provide a practical framework for optimal execution, allowing portfolios to be adjusted gradually over time while accounting for market impact and dynamic trading conditions~\cite{almgren2001optimal,obizhaeva2013optimal}. However, these advantages come at a computational cost. As the number of assets and trading periods increases, the resulting optimization problems become high-dimensional and computationally demanding. In practice, widely used CPU-based solvers can struggle to meet the latency requirements of modern trading environments, especially in multi-period settings where decisions must be updated continuously. This computational bottleneck limits the deployment of sophisticated portfolio optimization models in real-time applications.

In this technical report, we propose FlashFolio, a GPU-based solver for portfolio optimization. The key idea is to leverage the massive parallelism and high memory bandwidth of modern GPUs to accelerate the solution of large-scale convex optimization problems arising in both single-period and multi-period settings. Our approach is designed to exploit the structure of the problem through operations that map efficiently onto GPU architectures. This enables improvements in computational speed while maintaining solution accuracy, making it possible to solve large portfolio optimization problems at finer time scales in practice.

\section{Problem Formulation}
In this section, we introduce the portfolio optimization models studied in this work. We begin with a single-period formulation that determines an optimal post-trade portfolio from the current holdings by balancing expected return against risk and trading costs. We then extend this formulation to a multi-period setting, which models a dynamic execution process in which the portfolio is rebalanced gradually over time toward a prescribed terminal target.

\subsection{Single-period model}\label{sec:single-period}
We first consider a single-period portfolio optimization problem. Let $\mygmv$ denote the maximum gross market value (GMV) of the portfolio, which measures the total buying power available to the investor, and let $n$ denote the size of the investment universe, i.e., the number of tradable assets. The decision variable is $w \in \mathbb{R}^n$, where each entry of $w$ represents the fraction of GMV allocated to a corresponding asset. Given an initial portfolio $w_0 \in \mathbb{R}^n$, the goal is to determine a feasible post-trade portfolio that appropriately balances return, risk, and trading cost. We formulate the problem as
\begin{align}\label{eq:sppm obj}
        \max_{w}&\quad \mygmv\cdot\alpha^\top w - \lambda_1\mygmv^2\cdot w^\top \Sigma w - \lambda_2\mygmv\cdot \, s^\top |w - w_{0}| - \lambda_3\mygmv\cdot \, q^\top |w - w_{0}|^d\\
        \mathrm{s.t.}&\quad  L \leq A w \leq U, \label{eq:exposure}\\
        & \quad \|w\|_1 \le 1 \ . \label{eq:budget}
\end{align}

The objective \eqref{eq:sppm obj} consists of four terms corresponding to expected return, portfolio risk, linear transaction cost, and nonlinear market impact.
The first term models expected return, where $\alpha \in \mathbb{R}^n$ denotes the vector of predicted asset returns. The second term models portfolio risk through the quadratic form induced by the risk matrix $\Sigma \in \mathbb{R}^{n\times n}$. 
The third term represents linear transaction cost induced by bid--offer spreads, where $s \in \mathbb{R}_+^n$ is the vector of bid-offer spread coefficients. It penalizes the absolute trade vector $|w-w_0|$ and hence discourages excessive turnover. The fourth term represents nonlinear market impact, where $q \in \mathbb{R}_+^n$ is the vector of impact coefficients and $d \in (1,2]$ is the impact exponent. 
This term reflects the empirical fact that large trades in less liquid assets incur disproportionately higher execution costs. The nonnegative parameters $\lambda_1$, $\lambda_2$, and $\lambda_3$ are penalty coefficients that control the trade-off among risk, bid--offer cost, and market impact.

Constraint \eqref{eq:exposure} represents linear exposure constraints, where $A \in \mathbb{R}^{m\times n}$ is the exposure matrix, and $L,U \in \mathbb{R}^m$ are the corresponding lower and upper bounds. In practical portfolio construction, these constraints may encode market-neutrality requirements, sector or industry bounds, and limits on style-factor exposures. Constraint \eqref{eq:budget} is the total investment constraint, which ensures that the total amount invested, including both long and short positions, does not exceed the portfolio's available buying power.

Overall, the single-period model selects a feasible post-trade portfolio that balances expected return against multiple sources of risk and trading friction while satisfying investment and exposure limits.

\subsection{Multi-period model}\label{sec:multi-period}
The single-period model determines an optimal portfolio for a single rebalancing step. In many applications, however, trading is carried out over multiple rebalancing periods rather than in a single step. This is particularly important when transaction costs and market impact are significant, since it may be preferable to move gradually from the current portfolio to a target portfolio instead of executing the entire trade immediately.

To capture this dynamic setting, we extend the model to a multi-period formulation. Let $T$ denote the investment horizon, i.e., the total number of rebalancing periods. The initial portfolio $w_0 \in \mathbb{R}^n$ and the terminal portfolio $w_T \in \mathbb{R}^n$ are given, while the intermediate portfolios $w_1,\ldots,w_{T-1} \in \mathbb{R}^n$ are decision variables. Unlike the single-period model, which determines only one post-trade allocation, the multi-period model optimizes an entire rebalancing trajectory from $w_0$ to $w_T$. The multi-period portfolio optimization problem is formulated as

{\small\begin{align}
        \max_{w_1,...,w_{T-1}}&\quad \sum_{t=1}^{T} \mygmv\cdot\alpha_t^\top w_t - \lambda_1\mygmv^2\cdot (w_t - w_T)^\top \Sigma (w_t - w_T) - \lambda_2\mygmv\cdot \, s_t^\top |w_t - w_{t-1}| - \lambda_3\mygmv\cdot \, q_t^\top |w_t - w_{t-1}|^d \label{eq: mppm-obj}\\
        \mathrm{s.t.}&\quad  L \leq A w_t \leq U \quad \forall t=1, ...,T-1\ , \label{eq:mppm-expose}\\
        & \quad \|w_t\|_1 \le 1 \qquad \ \ \forall t=1, ...,T-1\ .\label{eq:mppm-budget}
\end{align}}
Here $\alpha_t$, $s_t$, and $q_t$ may vary across time to reflect changing return forecasts, bid--offer spreads, and liquidity conditions. The per-period objective retains the same basic ingredients as in the single-period model, but with an explicitly dynamic interpretation. The return term captures the benefit of holding portfolio $w_t$ during period $t$. The transaction-cost and market-impact terms depend on the trade increment $w_t-w_{t-1}$, reflecting the fact that trading cost is incurred when moving from one portfolio to the next. The quadratic term $(w_t-w_T)^\top \Sigma (w_t-w_T)$
penalizes deviations from the target portfolio along the execution path, thereby encouraging a gradual and risk-aware transition toward the desired terminal holdings.

The exposure constraints \eqref{eq:mppm-expose} and total investment constraints \eqref{eq:mppm-budget} are imposed at every intermediate period, ensuring that the portfolio remains feasible throughout the entire rebalancing process rather than only at the beginning and the end.

In summary, the multi-period model extends the single-period formulation by optimizing an entire rebalancing trajectory rather than a single post-trade allocation. This dynamic structure allows execution speed to adapt over time to the trade-off among return, risk, and trading costs, but it also significantly increases the computational complexity of the resulting optimization problem.

\section{Benchmark Dataset}\label{sec:data}

This section describes how we construct the benchmark instances used in our numerical experiments. Our goal is not to design a tradable forecasting or execution model, but to build realistic and controlled portfolio optimization instances from market data for evaluating solver performance. To this end, we specify how the alpha signal, risk model, and transaction cost coefficients are generated for the single-period model, and how these inputs are extended to the multi-period setting.

\subsection{Alpha Construction}

For each date $t$, we construct a synthetic alpha signal $\alpha \in \mathbb{R}^n$ from subsequent realized returns together with stock-specific noise. Specifically, for each stock $i$, we define
\[
\alpha_i
=
\frac{1}{5}\sum_{k=t}^{t+4} r_{k,i} + 0.5\,\epsilon_{t,i},
\]
where
\[
r_{k,i}
=
\frac{p_{k+1,i}-p_{k,i}}{p_{k,i}}
\]
denotes the daily return from date $k$ to $k+1$, and $p_{k,i}$ is the closing price of stock $i$ on date $k$. The perturbation term $\epsilon_{t,i}$ is drawn from a zero-mean Gaussian distribution with scale determined by the stock's daily volatility. Specifically, we use implied daily volatility when available, and otherwise replace it with the two-year historical standard deviation of daily returns.

This construction uses future returns only to generate a synthetic benchmark signal and is not intended to represent an implementable forecasting model. Its purpose is to produce realistic variation in the return term across benchmark instances while keeping the focus of the paper on optimization performance rather than alpha modeling.

The implied volatility data are obtained from the 30-day, 50-delta volatility curve in OptionMetrics via WRDS, and the stock price data are obtained from CRSP via WRDS.

\subsection{Risk Model Construction}

We construct the risk matrix using the factor-model representation
\[
\Sigma = \beta \Sigma_f \beta^\top + D,
\]
where $\beta \in \mathbb{R}^{n \times p}$ is the risk exposure matrix associated with $p$ underlying factors, $\Sigma_f \in \mathbb{R}^{p \times p}$ is the factor covariance matrix, and $D \in \mathbb{R}^{n \times n}$ is a diagonal matrix of specific variances. In our experiments, we use $p=6$ factors: the market factor together with the five Fama--French factors.

For each date $t$, we use the previous $T=504$ business days of factor returns to estimate both factor exposures and factor covariance. Let $f_k \in \mathbb{R}^p$ denote the factor return vector on day $k$, and let
\[
F =
\begin{bmatrix}
f_{t-T}^\top \\
\vdots \\
f_{t-1}^\top
\end{bmatrix}
\in \mathbb{R}^{T \times p}
\]
collect these factor returns over the estimation window. For each stock $i$, we estimate its factor loading vector $\beta_i \in \mathbb{R}^p$ by regressing its historical daily returns on the factor returns:
\[
r_i \sim F \beta_i,
\]
where $r_i \in \mathbb{R}^{T}$ denotes the vector of historical returns of stock $i$ over the same window. To improve robustness, factor exposures are estimated only for stocks with at least two months of return history.

The factor covariance matrix is estimated from the same window as
\[
\Sigma_f
=
\frac{1}{T-1}
\left(
F - \frac{1}{T}\mathbf{1}\mathbf{1}^\top F
\right)^\top
\left(
F - \frac{1}{T}\mathbf{1}\mathbf{1}^\top F
\right),
\]
where $\mathbf{1}\in\mathbb{R}^T$ is the vector of ones.

For each stock $i$, the residual return on day $k$ is then defined as
\[
\theta_{k,i} = r_{k,i} - \beta_i^\top f_k
\qquad k=t-T,\ldots,t-1,
\]
and the diagonal entries of the specific risk matrix are estimated by the sample variance of these residual returns:
\[
D_{ii}
=
\frac{1}{T-1}
\sum_{k=t-T}^{t-1}
\left(
\theta_{k,i} - \frac{1}{T}\sum_{j=t-T}^{t-1} \theta_{j,i}
\right)^2.
\]

All quantities in the risk model depend on the date $t$. The factor return data are obtained from the Fama-French dataset~\cite{FAMA19933}.

\subsection{Transaction Cost Construction}

We next construct the transaction cost coefficients used in the bid--offer spread and market impact terms.

The bid--offer cost vector $s$ is derived from bid and ask prices obtained from CRSP. The market impact coefficient vector $q$ is defined as
\[
q = \sigma \left(\frac{\mathrm{GMV}}{\mathrm{ADV}}\right)^{d-1},
\]
where $\sigma \in \mathbb{R}^n$ denotes daily volatility and $\mathrm{ADV} \in \mathbb{R}^n$ denotes average daily volume. For each stock $i$ at date $t$, the average daily volume is computed over the previous $60$ trading days as
\[
\mathrm{ADV}_{t,i}
=
\frac{1}{60}\sum_{k=t-60}^{t-1} \mathrm{DV}_{k,i},
\]
where $\mathrm{DV}_{k,i}$ is the daily trading volume of stock $i$ on day $k$.

For stocks without implied volatility data, we again replace $\sigma_{t,i}$ with the two-year historical standard deviation of daily returns. The daily trading volume data are obtained from CRSP, and the implied volatility data are obtained from the 30-day, 50-delta volatility curve in OptionMetrics.

\subsection{Multi-Period Setting}

We next describe how the single-period inputs are extended to construct multi-period benchmark instances.

For alpha, we assume a constant signal across periods and set
\[
\alpha_t = \alpha,
\qquad t=1,\ldots,T.
\]
For risk, we keep the covariance matrix fixed over the trading horizon and use the same matrix $Q$ at every period.

The transaction cost coefficients vary across periods according to a prescribed intraday liquidity profile. Specifically, we set
\[
s_t = \frac{\eta_t}{\sum_{k=1}^{T-1}\eta_k}\, s,
\qquad
q_t = \frac{\eta_t}{\sum_{k=1}^{T-1}\eta_k}\, q,
\]
where $\eta_t$ is an intraday scaling factor following a U-shaped pattern,
\begin{equation}\label{eq:ushape}
    \eta_t \coloneqq 0.5 + 2\left(\frac{t}{T}-0.5\right)^2.
\end{equation}
This profile assigns greater liquidity near the beginning and end of the horizon and lower liquidity in the middle, thereby mimicking the empirical U-shaped intraday pattern commonly observed in trading activity.

This construction keeps the alpha and risk inputs fixed across periods while allowing trading costs to vary over time. As a result, the multi-period benchmark isolates the computational effect of dynamic execution and time-coupled transaction costs, while remaining closely connected to the single-period instances.

\section{Numerical Experiments}\label{sec:compute}

In this section, we evaluate the performance of FlashFolio on portfolio optimization problems. We benchmark FlashFolio against MOSEK on both single-period and multi-period instances constructed from the dataset described in Section~\ref{sec:data}.

\subsection{Experiment setup}

\paragraph{Implementation.}
We implement FlashFolio in JAX 0.7.2 and compare its performance against the conic solver MOSEK 11.0.27~\cite{mosek}. For MOSEK, we construct the conic reformulation of the single-period and multi-period models in JuMP and provide the resulting models as input to MOSEK.

\paragraph{Computing Environment.}
FlashFolio is executed on an NVIDIA H100-SXM-80GB GPU with CUDA 12.4, while MOSEK is run on a Linux cluster equipped with dual Intel Xeon Platinum 8462Y+ CPUs, using 16 CPU cores and 128 GB of RAM. The time limit is set to 360 seconds for both solvers. FlashFolio is initialized at the initial portfolio $w_0$.

\paragraph{Benchmark dataset}
To evaluate the robustness of these solvers, we consider a broad range of risk-aversion parameters $\lambda_1 \in \{10^{-8}, 10^{-7}, 10^{-6}, 10^{-5}, 10^{-4}\}$ and transaction cost penalties $\lambda_2 = \lambda_3 \in \{1, 10, 100, 1000, 10000\}$ for both single-period model and multi-period model. For the multi-period model, we set the investment horizon to $T=10$ and use the U-shaped intraday profiles of average daily volume defined in \eqref{eq:ushape}, with ADV peaking at $t=0$ and $t=T$. 
For each parameter configuration, we generate 21 instances using monthly data from 2022-01-03 to 2023-08-31.

\textbf{Termination criteria.} 
FlashFolio terminates when its fixed-point residual falls below a prescribed tolerance $\epsilon$. In all experiments, we set $\epsilon = 10^{-8}$. MOSEK uses a different termination logic based on the KKT conditions of the conic reformulation, so the stopping criteria of the two solvers are not directly identical. For MOSEK, we use the default tolerances of $10^{-8}$ for primal feasibility, dual feasibility, and relative gap.

\paragraph{Evaluation metrics.} 
We evaluate solver performance based on the number of solved instances and the Shifted Geometric Mean (SGM) of solving time with a shift of $k=1$ seconds:
\[
\text{SGM1} = \left(\prod_{i=1}^n (t_i+1)\right)^{1/n}-1
\]
where $t_i$ is the solve time for the i-th instance.
We note that MOSEK returns the status \texttt{SLOW\_PROGRESS} on some instances. In practice, however, the returned solutions remain of sufficiently high quality, and we therefore include such instances in the solved count.

\subsection{Benchmark Results}

Table~\ref{tab:sgm1} compares the performance of FlashFolio and MOSEK on single-period and multi-period portfolio optimization instances under a range of parameter settings.

For the single-period model, both solvers successfully solve all 21 instances across all tested parameter configurations. In this regime, the main difference is therefore runtime rather than robustness. FlashFolio is uniformly faster than MOSEK, with speedups ranging from $2.4\times$ to $12.9\times$ depending on the parameter setting.

For the multi-period model, the advantage of FlashFolio becomes even more pronounced. FlashFolio solves all 21 instances in nearly all tested settings and maintains consistently low runtimes across the entire parameter range. By comparison, MOSEK begins to lose robustness in more challenging regimes, solving as few as 16 out of 21 instances in the hardest cases. At the same time, FlashFolio delivers substantially lower runtimes, with speedups ranging from $15.8\times$ to $48\times$.

Overall, these results show that FlashFolio provides not only a clear efficiency advantage over MOSEK, but also stronger numerical robustness on the more challenging multi-period instances.

\begin{table}[!htbp]
\centering
\caption{Performance comparison of FlashFolio and MOSEK on single-period and multi-period portfolio optimization instances. Count denotes the number of instances solved out of 21. Time denotes the shifted geometric mean of solve time in seconds with shift 1 (SGM1). Speedup is computed as the ratio of MOSEK time to FlashFolio time.}
\label{tab:sgm1}
\small
\setlength{\tabcolsep}{4pt}
\begin{tabular}{llrrrrrrrrrr}
\toprule
\multicolumn{2}{c}{Parameters} & \multicolumn{5}{c}{Single-period model} & \multicolumn{5}{c}{Multi-period model} \\
\cmidrule(lr){1-2} \cmidrule(lr){3-7} \cmidrule(lr){8-12}
\multicolumn{1}{l}{\multirow{2}{*}{$\lambda_2,\lambda_3$}} & \multicolumn{1}{l}{\multirow{2}{*}{$\lambda_1$}}
& \multicolumn{2}{c}{FlashFolio} & \multicolumn{2}{c}{MOSEK} & \multicolumn{1}{c}{\multirow{2}{*}{\textbf{Speedup}}}
& \multicolumn{2}{c}{FlashFolio} & \multicolumn{2}{c}{MOSEK} & \multicolumn{1}{c}{\multirow{2}{*}{\textbf{Speedup}}} \\
\cmidrule(lr){3-4} \cmidrule(lr){5-6} \cmidrule(lr){8-9} \cmidrule(lr){10-11}
& & Count & Time & Count & Time &  & Count & Time & Count & Time & \\
\midrule
\multirow{5}{*}{1}
& $10^{-4}$ & 21 & 0.24 & 21 & 0.71 & \textbf{3.0$\times$} & 21 & 0.88 & 21 & 13.88 & \textbf{15.8$\times$} \\
& $10^{-5}$ & 21 & 0.15 & 21 & 0.93 & \textbf{6.2$\times$} & 21 & 0.88 & 21 & 17.52 & \textbf{19.9$\times$} \\
& $10^{-6}$ & 21 & 0.14 & 21 & 1.80 & \textbf{12.9$\times$} & 21 & 0.93 & 21 & 21.84 & \textbf{23.5$\times$} \\
& $10^{-7}$ & 21 & 0.30 & 21 & 2.20 & \textbf{7.3$\times$} & 21 & 1.33 & 21 & 27.35 & \textbf{20.6$\times$} \\
& $10^{-8}$ & 21 & 0.64 & 21 & 1.77 & \textbf{2.8$\times$} & 21 & 1.77 & 21 & 30.87 & \textbf{17.4$\times$} \\
\cmidrule(lr){1-12}
\multirow{5}{*}{10}
& $10^{-4}$ & 21 & 0.25 & 21 & 0.79 & \textbf{3.2$\times$} & 21 & 0.96 & 21 & 18.79 & \textbf{19.6$\times$} \\
& $10^{-5}$ & 21 & 0.16 & 21 & 0.97 & \textbf{6.1$\times$} & 21 & 1.00 & 21 & 27.02 & \textbf{27.0$\times$} \\
& $10^{-6}$ & 21 & 0.16 & 21 & 1.75 & \textbf{10.9$\times$} & 21 & 1.06 & 21 & 24.28 & \textbf{22.9$\times$} \\
& $10^{-7}$ & 21 & 0.18 & 21 & 1.63 & \textbf{9.1$\times$} & 21 & 1.09 & 21 & 31.47 & \textbf{28.9$\times$} \\
& $10^{-8}$ & 21 & 0.18 & 21 & 1.67 & \textbf{9.3$\times$} & 21 & 1.08 & 21 & 29.56 & \textbf{27.4$\times$} \\
\cmidrule(lr){1-12}
\multirow{5}{*}{100}
& $10^{-4}$ & 21 & 0.30 & 21 & 0.94 & \textbf{3.1$\times$} & 21 & 1.06 & 21 & 29.25 & \textbf{27.6$\times$} \\
& $10^{-5}$ & 21 & 0.23 & 21 & 1.01 & \textbf{4.4$\times$} & 21 & 1.07 & 21 & 33.33 & \textbf{31.1$\times$} \\
& $10^{-6}$ & 21 & 0.22 & 21 & 1.22 & \textbf{5.5$\times$} & 21 & 1.10 & 21 & 31.22 & \textbf{28.4$\times$} \\
& $10^{-7}$ & 21 & 0.21 & 21 & 1.71 & \textbf{8.1$\times$} & 21 & 1.12 & 21 & 27.76 & \textbf{24.8$\times$} \\
& $10^{-8}$ & 21 & 0.22 & 21 & 1.71 & \textbf{7.8$\times$} & 21 & 1.11 & 21 & 40.93 & \textbf{36.9$\times$} \\
\cmidrule(lr){1-12}
\multirow{5}{*}{1000}
& $10^{-4}$ & 21 & 0.42 & 21 & 1.10 & \textbf{2.6$\times$} & 21 & 1.25 & 21 & 34.57 & \textbf{27.7$\times$} \\
& $10^{-5}$ & 21 & 0.43 & 21 & 1.16 & \textbf{2.7$\times$} & 21 & 1.26 & 21 & 27.34 & \textbf{21.7$\times$} \\
& $10^{-6}$ & 21 & 0.42 & 21 & 1.25 & \textbf{3.0$\times$} & 21 & 1.25 & 19 & 40.99 & \textbf{32.8$\times$} \\
& $10^{-7}$ & 21 & 0.44 & 21 & 1.27 & \textbf{2.9$\times$} & 21 & 1.26 & 17 & 45.65 & \textbf{36.2$\times$} \\
& $10^{-8}$ & 21 & 0.44 & 21 & 1.70 & \textbf{3.9$\times$} & 21 & 1.24 & 16 & 59.48 & \textbf{48.0$\times$} \\
\cmidrule(lr){1-12}
\multirow{5}{*}{10000}
& $10^{-4}$ & 21 & 0.63 & 21 & 1.49 & \textbf{2.4$\times$} & 21 & 1.44 & 21 & 28.47 & \textbf{19.8$\times$} \\
& $10^{-5}$ & 21 & 0.62 & 21 & 1.49 & \textbf{2.4$\times$} & 21 & 1.44 & 21 & 28.94 & \textbf{20.1$\times$} \\
& $10^{-6}$ & 21 & 0.61 & 21 & 1.69 & \textbf{2.8$\times$} & 21 & 1.44 & 20 & 34.96 & \textbf{24.3$\times$} \\
& $10^{-7}$ & 21 & 0.63 & 21 & 1.65 & \textbf{2.6$\times$} & 21 & 1.46 & 21 & 49.90 & \textbf{34.2$\times$} \\
& $10^{-8}$ & 21 & 0.63 & 21 & 1.66 & \textbf{2.6$\times$} & 21 & 1.45 & 17 & 33.73 & \textbf{23.3$\times$} \\
\bottomrule
\end{tabular}
\end{table}

\section{Conclusion}

We presented FlashFolio, a GPU-powered solver for single-period and multi-period portfolio optimization with factor-based risk modeling, bid--offer spread costs, and nonlinear market impact. Our computational results show that FlashFolio consistently outperforms MOSEK across a broad range of problem settings, achieving speedups in the single-period model and even larger gains in the multi-period model, while also maintaining stronger robustness on challenging instances. These results demonstrate that GPU-based optimization can help improve the practicality of large-scale portfolio optimization in latency-sensitive applications. 

\section*{Disclaimer}
The views and opinions expressed in this paper are solely those of the author and do not necessarily reflect the views, policies, or positions of Point72 Asset Management or its affiliates. This research was conducted in the author’s personal capacity and does not represent investment advice, trading strategies, or official research of Point72. No confidential or proprietary information of Point72 was used in the preparation of this paper.

\bibliographystyle{amsplain}
\bibliography{ref-papers}

\end{document}